\documentclass[11pt]{article}
\font\smallit=cmti12
\usepackage{amssymb,amsmath,latexsym,amsthm}

\newenvironment{frcseries}{\fontfamily{frc}\selectfont}{}

\newcommand\md[1]{\,\,(\mathrm{mod }\, #1)}
\def\w{w_{\mathrm{\mathfrak{z}}}}
\def\wtwo{w_{\mathrm{\mathfrak{z,2}}}}
\def\wm{w_{\mathrm{\mathfrak{m}}}}
\def\wmtwo{w_{\mathrm{\mathfrak{m,2}}}}
\def\B{\hfill $\Box$ \vskip 5pt}

\newtheorem{thm}{Theorem}

\newtheorem{prop}[thm]{Proposition}
\theoremstyle{definition}
\newtheorem{defn}[thm]{Definition}

\makeatletter

\renewcommand\section{\@startsection {section}{1}{\z@}%
{-30pt \@plus -1ex \@minus -.2ex}%
{2.3ex \@plus.2ex}%
{\normalfont\normalsize\bfseries}}

\renewcommand\subsection{\@startsection{subsection}{2}{\z@}%
{-3.25ex\@plus -1ex \@minus -.2ex}%
{1.5ex \@plus .2ex}%
{\normalfont\normalsize\bfseries}}
\renewcommand{\@seccntformat}[1]{\csname the#1\endcsname. } %\quad}

\makeatother

\begin{document}

\begin{center}
\uppercase{\bf  Zero-sum Analogues of van der Waerden's Theorem on Arithmetic Progressions}\\[25pt]
{\bf Aaron Robertson} \\
{\smallit Department of Mathematics,
Colgate University, Hamilton, New York}\\
{\tt arobertson@colgate.edu}
\end{center}

\vskip 10pt
\centerline{\bf Abstract}
\noindent
Let $r$ and $k$ be  positive integers
with $r \mid k$.  Denote by $\w(k;r)$ the minimum integer such
that every coloring $\chi:[1,\w(k;r)] \rightarrow \{0,1,\dots,r-1\}$ admits a $k$-term arithmetic progression $a,a+d,\dots,a+(k-1)d$
with $\sum_{j=0}^{k-1} \chi(a+jd) \equiv 0 \md{r}$.  We investigate these numbers
as well as a ``mixed" monochromatic/zero-sum analogue.  We also present an interesting
reciprocity between the van der Waerden numbers and $\w(k;r)$.

\baselineskip=14pt

\section{Introduction}
Van der Waerden's theorem \cite{vdw} on arithmetic progressions states that
for $k,r \in \mathbb{Z}^+$, there exists a minimum integer $w(k;r)$   such that {every} $r$-coloring of $[1,w(k;r)]$ admits
a monochromatic $k$-term arithmetic progression.  
The determination of these numbers is notoriously difficult;
in fact, only seven of these numbers are known.

In this article we investigate some zero-sum analogues of van der Waerden's
theorem.

\begin{defn} Let $a_1,a_2,\dots,a_n$ be a sequence of non-negative integers
and let $m \in \mathbb{Z}^+$.  We say that the sequence is {\it $m$-zero-sum} if
$\sum_{i=1}^n a_i \equiv 0 \md{m}$.
\end{defn}

The seminal result in the area of zero-sum sequences is the
Erd\H{o}s-Ginzberg-Ziv theorem \cite{EGZ}, which states that
any sequence of $2n-1$ integers contains an $n$-zero-sum subsequence of
$n$ integers.  Since around 1990, research activity concerning zero-sum
results has flourished,   through both the lens of additive number theory and
Ramsey theory.  For example,
the weighted Erd\H{o}s-Ginzberg-Ziv theorem due to Grynkiewicz \cite{G} allows us to multiply the integers
in the Erd\H{o}s-Ginzberg-Ziv theorem by weights.  This result states,
in particular, that if $w_1,w_2,\dots,w_n$ is an $n$-zero-sum sequence and
$a_1,a_2,\dots,a_{2n-1}$ is a sequence of $2n-1$ integers, then there exists
an $n$-term subsequence  $a_{i_1}, a_{i_2},\dots,a_{i_n}$ and a permutation
$\pi$ of $\{i_1,i_2,\dots,i_n\}$ such that
$\sum_{j=1}^n w_j a_{\pi(i_{j})} \equiv 0\md{n}$.
Further recent results can be found in  \cite{AM}, \cite{BCRY}, and \cite{GG}
among many others.

Most investigations of zero-sum sequences do not have a structure imposed on them.  This is in contrast to zero-sum results on edgewise colored graphs,
which have been around for many years (see, e.g., \cite{AC}, \cite{B}, \cite{BD}, and \cite{C} ).
Some notable exceptions are found in works of Bialostocki, such as \cite{BBS} and \cite{BSS} where the zero-sum sequence $x_1,x_2,\dots,x_n$ satisfies
$\sum_{i=1}^{n-1} x_i < x_n$ and in \cite{BBCY} where $x_{i+1}-x_{i} \leq x_i-x_{i-1}$ for $1 \leq i \leq n-1$.
These exceptions, however, do not have a rigid structure imposed on them due to the use of
inequality.
In this article we investigate zero-sum arithmetic progressions, thereby imposing a rigid structure on the sequences.  In  subsequent article \cite{R} we investigate zero-sum sequences with a different rigid structure (where one term is the sum of all other terms).

The first two zero-sum analogues of van der Waerden's
theorem we will investigate are given in the following definitions.

\begin{defn} Let $k$ and $r$ be  positive integers
such that $r \mid k$.  We denote by $\w(k;r)$ the minimum integer such
that every coloring $[1,\w(k;r)]$ with the colors $0,1,\dots,r-1$ (which we
may refer to as $\mathbb{Z}_r$) admits a $k$-term $r$-zero-sum arithmetic progression.
\end{defn}

\begin{defn} Let $k$ and $r$ be positive integers
such that $r \mid k$.  We denote by $\wtwo(k;r)$ the minimum integer such
that every coloring $[1,\wtwo(k;r)]$ with the colors $0$ and $1$ admits a $k$-term $r$-zero-sum arithmetic progression.
\end{defn}

Implicit in the above definitions is the existence of the respective minimum numbers, which follows directly from the existence of $w(k;r)$.  Note that
we need only prove the existence of $\w(k;r)$ since we easily
have $\wtwo(k;r) \leq \w(k;r)$ as $\mathbb{Z}_2 \subseteq \mathbb{Z}_r$.  The existence of $\w(k;r)$ comes
from $\w(k;r) \leq w(k;r)$ as any $k$-term monochromatic arithmetic
progression is $r$-zero-sum when $r \mid k$.  When $r \nmid k$, 
coloring every integer of $\mathbb{Z}^+$ with the color $1$ does not
admit a $k$-term $r$-zero-sum arithmetic progression.
In an interesting turn of events, we will see later in this article that the independent
existence of $\w(k;r)$ implies the existence of $w(k;r)$.

\section{Some Computation}

We start with results from computer calculations.  We wrote the Fortran programs
{\tt ZSAP.f} and {\tt ZSAP2.f}, available at
{\tt http://www.aaronrobertson.org}, for the determinations of $\w(k;r)$ and
$\wtwo(k;r)$, respectively, for small values of $k$ and $r$.
The algorithm used in both is a standard backtrack model to exhaustively search
the colorings for $k$-term $r$-zero-sum arithmetic progressions.

Based on the values (see Tables 1 and 2, below) we see some interesting patterns.
For $r=2$, we clearly have $\w(k;2)=\wtwo(k;2)$ by definition, but it appears
that $\w(k;2)=2k-1$.  We prove this in the next section.  (It is interesting to note
that $2k-1$ is the formula from the Erd\H{o}s-Ginzberg-Ziv zero-sum theorem.)
For $k=3,6,9,12$, we have $\w(k;3)=\wtwo(k;3)=k^2$ and we investigate this in the
next section as well.  Along the diagonal, we see familiar van der Waerden numbers
appear.  This is explained in Section \ref{sec3} as well.

\begin{center}
\begin{tabular}{||l|c|c|c|c||} \hline\hline
\multicolumn{1}{||c|}{${}_{\mbox{\large $k$}} \,\,\,\,{}_\diagdown {\,\,\mbox{\large $r$}}$}      &       \multicolumn{1}{c|}{2}
   &  \multicolumn{1}{c|}{3}  &  \multicolumn{1}{c|}{4}  &
    \multicolumn{1}{c||}{5}  \\[5pt] \hline
    $2$                   &  3  &  $\infty$  &  $\infty$  &  $\infty$      \\[-1pt] 
    $3$                   &  $\infty$   &  $9$  &  $\infty$  &  $\infty$        \\[-1pt] 
    $4$                   &  7  &  $\infty$  &  $35$  &  $\infty$    \\[-1pt] 
    $5$                   &  $\infty$  &  $\infty$  &   $\infty$ &$\geq 294$     \\[-1pt] 
    $6$                   &  11  &  $36$  &  $\infty$  &  $\infty$      \\[-1pt] 
    $7$                   &  $\infty$  &  $\infty$  &  $\infty$  &  $\infty$    \\[-1pt] 
    $8$                   &  15  &  $\infty$  &  $\geq 108$  &  $\infty$     \\[-1pt] 
    $9$                   &  $\infty$  &  $81$  &  $\infty$  &  $\infty$       \\[-1pt] 
    $10$                  &  19  &  $\infty$  &  $\infty$  &  $\geq ???$      \\[-1pt] 
    $11$                  &  $\infty$  &  $\infty$  &  $\infty$  &  $\infty$      \\[-1pt] 
    $12$                  &  23  &  $144$  &  $\geq ??$  &  $\infty$     \\
\hline\hline
    \end{tabular}
\vskip 5pt
\centerline{\small {\bf Table 1}:  Values and lower bounds for $\w(k;r)$ for small $k$ and $r$}
    \end{center}

\noindent {\bf Remark.} The lower bounds for $\w(8;4)$ and
$\w(5;5)$ were obtained within a few hours but were not improved
upon after 770 hours of computation.  The lower bounds for
$\w(12;4)$ and $w(10;5)$ were obtained after ???? hours of computation,
with no extended time used to try to improve upon them.

\begin{center}
\begin{tabular}{||l|c|c|c|c||} \hline\hline
\multicolumn{1}{||c|}{${}_{\mbox{\large $k$}} \,\,\,\,{}_\diagdown {\,\,\mbox{\large $r$}}$}      &       \multicolumn{1}{c|}{2}
   &  \multicolumn{1}{c|}{3}  &  \multicolumn{1}{c|}{4}  &
    \multicolumn{1}{c||}{5}  \\[5pt] \hline
    $2$                   &  3  &  $\infty$  &  $\infty$  &  $\infty$      \\ 
    $3$                   &  $\infty$   &  $9$  &  $\infty$  &  $\infty$       \\
    $4$                   &  7  &  $\infty$  &  $35$  &  $\infty$    \\
    $5$                   &  $\infty$  &  $\infty$  &   $\infty$ &$178$      \\
    $6$                   &  11  &  $36$  &  $\infty$  &  $\infty$     \\
    $7$                   &  $\infty$  &  $\infty$  &  $\infty$  &  $\infty$   \\
    $8$                   &  15  &  $\infty$  &  $80$  &  $\infty$     \\
    $9$                   &  $\infty$  &  $81$  &  $\infty$  &  $\infty$       \\
    $10$                  &  19  &  $\infty$  &  $\infty$  &  $\geq ???$     \\
    $11$                  &  $\infty$  &  $\infty$  &  $\infty$  &  $\infty$     \\
    $12$                  &  23  &  $144$  &  $244$  &  $\infty$    \\
\hline\hline
    \end{tabular}
\vskip 5pt
\centerline{\small {\bf Table 2}:  Values and lower bounds for $\wtwo(k;r)$ for small $k$ and $r$}
    \end{center}

\noindent
{\bf 
Remark.} The determination of $\wtwo(12;4)$ took about 8 days, while the lower bound
for $\wtwo(10;5)$ was established after only a few minutes of searching,
but further impovement was not achieved after ?? hours.
All other  values took less than a few hours.

\section{Formulas}\label{sec3}

As previously mentioned, some interesting patterns can be seen in Tables 1 and 2.
In this section we explore these patterns.  We start with a formula for the $r=2$ columns
and note that it also proves the existence of $\w(k;2)$ for all even $k$, independent
of the existence of the van der Waerden number $w(k;2)$.

\begin{thm} Let $k \in \mathbb{Z}^+$ be even.  Then $\w(k;2)=\wtwo(k;2)=2k-1$.
\end{thm}

\noindent
{\it Proof.} The first equality is by definition.  To show that $\w(k;2) = 2k-1$ we will provide
matching upper and lower bounds.  First, it is easy to check that the $2$-coloring of $[1,2k-2]$ with all
integers colored 0 except for integer $k$ avoids $k$-term $2$-zero-sum arithmetic progressions
as any such arithmetic progression must consist of $k$ consecutive integers and,
hence, exactly one integer of color $1$.  Hence, $\w(k;2) \geq 2k-1$.

We next show that $\w(k;2) \leq 2k-1$ by contradiction, assuming that there exists a
coloring $\chi$ of $[1,2k-1]$ by $\mathbb{Z}_2$ with no $k$-term $2$-zero-sum arithmetic progression.
Let $A=\{1,3,5,\dots,k-1\}$ and $B=\{k+1, k+3,\dots,2k-1\}$.  Since $A\cup B$ is a $k$-term
arithmetic progression, we assume that the sum of the colors of the integers in $A \cup B$
is odd.  Hence, one of $A$ and $B$ has an even number of integers of color 1, while the other
has an odd number of integers of color 1.  Without loss of generality, let
$A$ have an even number of integers of color 1.

Consider $S(x) = \sum_{i=x}^{k+x-1} \chi(i)$ for $x \in [1,k]$.  Next note that
$S(x+1)-S(x) = \chi(k+x) - \chi(x)$ for $x \in [1,k-1]$.  Since we assume that $S(x) \equiv 1 \md{2}$ for all $x \in [1,k]$
we must have $\chi(k+x) = \chi(x)$ for $x \in [1,k-1]$.  This means that $[k+1,2k-1]$ is colored in the exact
same way as $[1,k-1]$.  This contradicts our determination that $A$ has an even number of integers
of color 1, while $B$ has an odd number of integers of color 1. \B

For the $r=3$ column, we can justify lower bounds that match all of the calculated numbers
via  Theorems \ref{th5} and \ref{th6}, below.

\begin{thm} \label{th5}
 Let $k \in \mathbb{Z}^+$ with $3 \mid k$. If
$k+1$ is prime, then $\wtwo(k;3)\geq k^2$.
\end{thm}

\noindent
{\it Proof.}  We will show that the $2$-coloring $\chi$ of $[1,k^2-1]$ defined
by 
$$(0^{k-1}11)^{k-1}$$ (i.e., the color pattern of $k-1$ consecutive
0s followed by two 1s, repeated $k-1$ times) avoids $k$-term
$3$-zero-sum arithmetic progressions.  

Consider an arbitrary
$k$-term arithmetic progression $a, a+d, a+2d,\dots,a+(k-1)d$.
Note that $a+(k-1)d \leq k^2-1$ gives us that $d \leq k$.
Since $k+1$ is prime, we have $(d,k+1)=1$.  It follows that
$\{a,a+d,\dots,a+(k-1)d\}$ when reduced modulo $k+1$ is a set of $k$ (distinct) elements of $\mathbb{Z}_{k+1}$.

Looking at our coloring, we interpret it as
$\chi(x) = 1$ if $x \equiv 0 \mbox{ or } k \md{k+1}$ and $\chi(x)=0$
otherwise.  Since our arithmetic progression hits $k$ distinct residues
modulo $k+1$ we see that $\sum_{j=0}^{k-1} \chi(a+jd) = 1 \mbox{ or } 2$
so that it is not $3$-zero-sum.
\B

\begin{thm} \label{th6}
Let $k \in \mathbb{Z}^+$ with $3 \mid k$.  
If $k+1=2p$ with $p$ prime, then $\wtwo(k;3)\geq k^2$.
\end{thm}

\noindent
{\it Proof.}  We will show that the $2$-coloring $\chi$ of $[1,k^2-1]$ defined
by $$(0^{p-2}101^{p-2}01)^{k-1}$$ avoids $k$-term
$3$-zero-sum arithmetic progressions.  

Consider an arbitrary
$k$-term arithmetic progression $a, a+d,\dots,a+(k-1)d$.
Note that $a+(k-1)d \leq k^2-1$ gives us that $d \leq k$.
Hence, $(d,k+1) \in \{1, 2, p\}$.  We will determine
$\sum_{j=0}^{k-1} \chi(a+jd) \md{3}$ based on
the value of $(d,k+1)$.  We will also use the
fact that $p \equiv 2 \md{3}$ which follows from
$k=2p-1 \equiv 0 \md{3}$.

\vskip 5pt
\noindent
{\tt Case 1.} $(d,k+1)=1$.
It follows that
$\{a,a+d,\dots,a+(k-1)d\}$ when reduced modulo $k+1$ is a set of $k$ 
(distinct) residues of $\mathbb{Z}_{k+1}$.  Hence, $\sum_{j=0}^{k-1} \chi(a+jd) =
p-1 \mbox{ or } p$.  Since $p \equiv 2 \md{3}$ we have
$\sum_{j=0}^{k-1} \chi(a+jd) \equiv 1 \mbox{ or } 2 \md{3}$
so that our arithmetic progression is not $3$-zero-sum.
\hfill $\diamond$

\vskip 5pt
\noindent
{\tt Case 2.} $(d,k+1)=2$.  By reducing all terms of $\{a,a+d,\dots,a+(k-1)d\}$ modulo $k+1$ we see that we have either $\{0,2,4,\dots,k-1\}$ or
$\{1,3,5,\dots,k\}$.  Looking at $0^{p-2}101^{p-2}01$ we see that the coloring
of the even terms is $0^{\frac{p-3}{2}}
1^{\frac{p+3}{2}}$,
while the coloring of  the odd terms is $0^{\frac{p+1}{2}} 1^{\frac{p-3}{2}}0$. 

Next, when reducing all terms of $\{a,a+d,\dots,a+(k-1)d\}$ modulo $k+1$ we
see that every residue except for one  of $\{0,2,4,\dots,k-1\}$
is congruent modulo $k+1$ to precisely two terms of the arithmetic
progression.  The same holds for the residue set $\{1,3,5,\dots,k\}$.

In the situation where $\{a,a+d,\dots,a+(k-1)d\}$ modulo $k+1$ is $\{0,2,4,\dots,k-1\}$ we have
$\sum_{j=0}^{k-1} \chi(a+jd) = 2(\frac{p+3}{2}) + \epsilon$ where 
$\epsilon \in \{-1,0\}$.  We have $2(\frac{p+3}{2}) + \epsilon
\equiv p + \epsilon \md{3}$.  Since $p \equiv 2 \md{3}$, we see that
our arithmetic progression is not $3$-zero-sum in this situation.

In the situation where $\{a,a+d,\dots,a+(k-1)d\}$ modulo $k+1$ is $\{1,3,5,\dots,k\}$ we have\break
$\sum_{j=0}^{k-1} \chi(a+jd) = 2(\frac{p-3}{2}) + \epsilon$ where 
$\epsilon \in \{-1,0\}$.  We have $2(\frac{p-3}{2}) + \epsilon
\equiv p + \epsilon \md{3}$.  Since $p \equiv 2 \md{3}$, we see that
our arithmetic progression is not $3$-zero-sum in this situation.
\hfill $\diamond$

\vskip 5pt
\noindent
{\tt Case 3.} $(d,k+1)=p$.  We must have $d=p$ in this situation since
$d \leq k$.  Looking at $0^{p-2}101^{p-2}01$ as a coloring of
$[1,k+1]=[1,2p]$ we see that $\chi\left(p+i\right) = \chi(i) + 1 \md{2}$
for $i=1,2,\dots,p$.  In particular,
$\chi(x) + \chi(x+p) = 1$ for any $x \in \left[1,p\right]$.
Given that our coloring $0^{p-2}101^{p-2}01$ is repeated $k-1$ times,
we have $\chi(x) + \chi(x+p) = 1$ for any $x$
where, for $\bar{x} \equiv x \md{k+1}$, we have $\bar{x}\in \left[1,p\right]$.

If $a \md{k+1}$ is between $1$ and $p$, inclusive, then we have
$$
\sum_{j=0}^{k-1} \chi(a+jd) = \sum_{j=0}^{k-1} \chi(a+jp)
= \chi(a+(k-1)p) + \sum_{\stackrel{0 \leq j \leq k-3}{j \mathrm{ \,\,even}}} \big(\,\chi(a+jp) + \chi(a+(j+1)p)\,\big) 
$$
$$
\hspace*{14pt}=  \chi(a+(k-1)p) + 1 \cdot \frac{k-1}{2}.
$$

$$ \hspace*{-5pt}
= \chi(a+(k-1)p) + p -1.
$$

Now, since $p \equiv 2\md{3}$, we see that regardless of
the value of $\chi(a+(k-1)p)$ we have
$\sum_{j=0}^{k-1} \chi(a+jd) \not \equiv 0 \md{3}$ so that
our arithmetic progression is not $3$-zero-sum.

If $a \md{k+1}$ is between $p+1$ and $k+1$, inclusive, then we have
$$
\sum_{j=0}^{k-1} \chi(a+jd) = \sum_{j=0}^{k-1} \chi(a+jp)
= \chi(a)+\sum_{\stackrel{1 \leq j \leq k-2}{j \mathrm{ \,\,odd}}} \big(\chi(a+jp) + \chi(a+(j+1)p)\big) 
$$
$$
\hspace*{13pt}= \chi(a) + 1 \cdot \frac{k-1}{2} .
$$

$$ \hspace*{-3pt}
= \chi(a) + p-1 .
$$

Again, since $p \equiv 2\md{3}$,   regardless of
the value of $\chi(a)$ 
our arithmetic progression is not $3$-zero-sum.
\hfill $\diamond$

As we have exhausted all possible values of $(d,k+1)$ and
shown that no $k$-term
$3$-zero-sum arithmetic progression exists under our
coloring in each situation, we are done.
\B

\vskip 5pt
\noindent
{\bf Remark.} Theorems \ref{th5} and \ref{th6} also give lower
bounds for $\w(k;3)$ with $k$ being a prime or twice a prime.
\vskip 5pt

The next proposition explains the appearance of the van der Waerden numbers along the
main diagonal of Table 2.  These same numbers are lower bounds for the main
diagonal of Table 1, where we see divergence occurring at $\w(5;5)$.

\begin{prop} \label{exist} Let $k \in \mathbb{Z}^+$.  Then $w(k;2)=\wtwo(k;k) \leq \w(k;k)$, 
\end{prop}

\noindent
{\it Proof.} The equality $w(k;2)=\wtwo(k;k)$ follows from the fact that the only way
a $k$-term sequence of $0$s and $1$s can be $k$-zero-sum is if all terms are 0s or all
terms are 1s, i.e., monochromatic.
Since $\mathbb{Z}_2 \subseteq \mathbb{Z}_r$, we have $\wtwo(k;r) \leq \w(k;r)$ so we are done.
\B

\section{A ``Mixed" Monochromatic/Zero-sum Analogue}

In this section we investigate an interplay between monochromatic and zero-sum arithmetic progressions.
We start with the question of whether or not by avoiding certain monochromatic arithmetic
progressions we can guarantee certain zero-sum arithmetic progressions.  To this end, consider the
following definition.

\begin{defn}
Let $k,\ell,r \in \mathbb{Z}^+$  with $k,\ell,r \geq 2$.  Define $\wm(k,\ell;r)$ to be the minimum integer $n$ such that
any coloring of $[1,n]$ by $\mathbb{Z}_r$ admits either a $k$-term monochromatic arithmetic
progression of a color other than $0$ or an $\ell$-term $r$-zero-sum arithmetic progression.
\end{defn}

Inherent in this definition is the existence of $\wm(k,\ell;m)$ for all positive integers $k,\ell,$ and $m$, so
we must justify this existence.  The existence follows easily from van der Waerden's theorem.
For the situation when $k \geq \ell$, we know that any coloring of $[1,w(k;r)]$ admits a monochromatic
$k$-term arithmetic progression.  If the color is anything but color $0$, then we are done, so we
assume that it is of color $0$.  But then we have an $\ell$-term arithmetic progression of color $0$,
which is necessarily $r$-zero-sum for any $r$.  For the situation when $k < \ell$, any coloring
of $[1,w(\ell;r)]$ admits a monochromatic $\ell$-term arithmetic progression.  If this color is $0$, then
it is $r$-zero-sum.  Otherwise, it contains a $k$-term monochromatic arithmetic progression of color
other than $0$.

Having the existence of $\wm(k,\ell;r)$, we see that these ``mixed" monochromatic/zero-sum
numbers, in particular $\wm(k,k;r)$, address the non-existence of $\w(k;r)$ when $r \nmid k$ (recall that the counterexample
was coloring all integers with color $1$).

Using the Fortran programs {\tt MZSAP.f} and {\tt MZSAP2.f}, available at {\tt www.aaronrobertson.org},
we have calculated the values in Tables 3 and 4.

\begin{center}
\begin{tabular}{l|c||c|c|c|c||} 
\multicolumn{1}{c}{}&\multicolumn{1}{c||}{}      &       
\multicolumn{1}{c|}{$\ell=2$}
   &  \multicolumn{1}{c|}{$\ell=3$}  &  \multicolumn{1}{c|}{$\ell=4$}  &
    \multicolumn{1}{c||}{$\ell=5$}  \\ \hline\hline  
\multicolumn{1}{c|}{$k$}&\multicolumn{1}{c||}{$r$}&\multicolumn{1}{c}{}&\multicolumn{1}{c}{}&\multicolumn{1}{c}{}&\\\hline
    & $2$                   &  3  &  6  &  7  &  10      \\ 
   2& $3$                  & 3 & 7 & 7 & 15        \\ 
    &$4$                   & 4 & 8 & 12& 20       \\\hline 
    & $2$                  & 3 & 7& 7&15        \\ 
   3& $3$                  & 7 &9 &16 & 25       \\ 
    &$4$                   & 7 & 21&28 &  47      \\\hline 
    &$2$                   & 3 & 8&  7& 20       \\ 
   4&$3$                   & 7 & 9 &18 & 33       \\ 
    &$4$                   & 11 &53 & 35& $\geq 97$       \\\hline 
    &$2$                   & 3 & 8& 7 & 21      \\ 
   5&$3$                   & 10 &9 &21 & 33       \\ 
    &$4$                   & 15 &219 &$35$ & $\geq ???$        \\ 
\hline\hline
    \end{tabular}
\vskip 5pt
\centerline{\small {\bf Table 3}:  Values and lower bounds for $\wm(k,\ell;r)$ for small $k, \ell$, and $r$}
    \end{center}

\noindent
{\bf Remark.} All exact values were achieved within a few hours of
computation time.  The lower bound for $\wm(4,5;4)$ was attained
quite quickly but not improved upon after 370 hours of computation.
The lower bound for $\wm(5,5;4)$ was reached after ??? hours, with
no computational  effort to improve upon it.

Examining Table 3, patterns do not pop out as they did in Tables 1 and 2.
There seems to be different behavior for a given $k$ depending on, perhaps,
the value of the gcd$(\ell, r)$ (see, e.g., the rows for $k=4,5$).  
We do see that for $k=2,3,$ and $4$ we have $\wm(k,k;k)=w(k;2)$, and
for $\ell=3,4$, and $5$ we have $\wm(3,\ell)=\ell^2$.  However, further calculation
shows that $\wm(3,6) = 33 \neq 6^2$.  
We can, however, provide formulas for the first three rows.

\begin{thm} \label{th9} Let $\ell \geq 2$ be an integer. Then 
$$
\wm(2,\ell;2)=\left\{
\begin{array}{ll}
2\ell -1&\mbox{ if } 2\mid \ell\\
2\ell&\mbox{ if }2 \nmid \ell
\end{array};
\right.
$$
$$
\wm(2,\ell;3)=\left\{
\begin{array}{ll}
2\ell -1&\mbox{ if } \ell \equiv 0,2,4 \md{6}\\
3\ell-2&\mbox{ if } \ell \equiv 3 \md{6}\\
3\ell&\mbox{ if } \ell \equiv 1,5 \md{6}
\end{array};
\right.
$$
and
$$
\wm(2,\ell;4) = \left\{
\begin{array}{ll}
3\ell&\mbox{ if } \ell \equiv 0,2,3,4 \md{6}\\
4 \ell&\mbox{ if } \ell \equiv 1,5 \md{6}
\end{array}.
\right.
$$
\end{thm}

\noindent
{\it Proof.} {\boldmath $\wm(2,\ell;2)$}. We first consider $\wm(2,\ell;2)$.  If $\ell$ is even, consider the coloring $0^{\ell -1}10^{\ell-2}$;
if $\ell$ is odd, consider the coloring $0^{\ell -1}10^{\ell-1}$.  It is routine to check that these do not admit
$2$-term arithmetic progressions of color $1$ or $\ell$-term $2$-zero-sum arithmetic progressions.  Now, 
with $\ell$ even, let
$n = 2\ell-1$ and assume, for a contradiction, that there exists a $2$-coloring of $[1,n]$ that avoids the
requisite arithmetic progressions.  Clearly, we may have at most one integer of color 1.  Furthermore, we can
have at most $\ell-1$ consecutive integers of color 0.  Hence, the only way to have a valid $2$-coloring of $[1,n]$
is with $0^{\ell-1}10^{\ell-1}$.  But then $1,3,5,\dots,2\ell-1$ is an $\ell$-term $2$-zero-sum arithmetic
progression, a contradiction.  The case when $\ell$ is odd is easier since, by the same reasoning, we can
only have $2\ell-1$ integers colored (by $0^{\ell-1}10^{\ell-1}$) and avoid the relevant arithmetic progression.
Hence,  every $2$-coloring of $[1,2\ell]$   admits either a
$2$-term arithmetic progressions of color $1$ or an $\ell$-term $2$-zero-sum arithmetic progressions.

\vskip 5pt
\noindent
{\boldmath $\wm(2,\ell;3)$}. Next, consider $\wm(2,\ell;3)$ and let $\ell
\equiv 0,2,4 \md{6}$ so that $\ell$ is even.  Again, it is easy to see that $0^{\ell -1}10^{\ell-2}$
does not admit $2$-term arithmetic progressions of color $1$ or $2$, or $\ell$-term $2$-zero-sum arithmetic progressions.
Assume, for a contradiction, that $\chi:[1,2\ell-1] \rightarrow \{0,1,2\}$ does not admit the relevant
arithmetic progressions.  Then at most one integer has color $1$ and at most one integer has color $2$.
Hence, if we use both colors $1$ and $2$ then our coloring has form $0^s10^t20^u$ or its reverse.  The following
argument works for either form, so we will assume we have $0^s10^t20^u$.  We know that $s,t,u \leq \ell-1$, subject to $s+t+u = 2\ell-3$.  
However, if $t \leq \ell-2$ then $0^x10^t20^y$ with $x+y+t=\ell-2$ is an $\ell$-term $3$-zero-sum arithmetic progression.
Since $s+t+u = 2\ell-3$ while $10^t2$ is at most $\ell$ terms, we see that $s+u \geq \ell - 1$, which
gives us that for some $x,y \geq 0$ we indeed have 
the existence of $0^x10^t20^y$
with $x+y+t=\ell-2$
contained in $\chi$, a contradiction.
If $\chi$ does not use both colors $1$ and $2$, then our coloring has form $0^s10^t$ or $0^s20^t$.  We will assume the former (the argument is the same for the latter).  We must have $s+t=2\ell-2$ and, hence, $s=t=\ell-1$.  But then
$1,3,5,\dots,2\ell-1$ is an $\ell$-term $3$-zero-sum arithmetic progression, a contradiction.
Hence, $\wm(2,\ell;3)=2\ell-1$ for $\ell$ even.

We next look at $\wm(2,\ell;3)$ for $\ell \equiv 3 \md{6}$.  For a lower bound, consider the coloring of $[1,3\ell-3]$ given by $0^{\ell-1}10^{\ell-1}20^{\ell-3}$.
Clearly we have no $2$-term arithmetic progression of color 1 or 2 and we do not have $\ell$ consecutive
integers that are $3$-zero-sum.  Hence, the only possible $\ell$-term $3$-zero-sum arithmetic progression must
have common gap 2.  Since $\ell \equiv 3 \md{6}$ we see that $\ell$ is odd.  This means that any arithmetic
progression with common gap 2 cannot contain both the color 1 and color 2.  Let $a,a+2,\dots,a+2(\ell-1)$ be
any arbitrary arithmetic progresion with common gap 2.
In order to have $a+2(\ell-1) \leq 3(\ell-1)$ we must have $a \leq \ell-1$.  This means that one of the
terms must have either color 1 or color 2, but that both colors 1 and 2 cannot occur in the progression.
Hence, $\sum_{i=0}^{\ell-1} \chi(a+2i) = 1$ or $2$ and is not $3$-zero-sum.  We conclude that $\wm(2,\ell;3) \geq 3\ell-2$.

To show that $\wm(2,\ell;3) \leq 3\ell-2$ for $\ell \equiv 3 \md{6}$, assume that $\chi:[1,3\ell-2] \rightarrow \{0,1,2\}$ is a coloring with no
$2$-term arithmetic progression of color 1 or 2 and no $\ell$-term $3$-zero-sum arithmetic progression.
We easily see that $\chi$ must use all 3 colors for otherwise we cannot have more than $2\ell-1$ integers of colors
only 0 and 1 (or 0 and 2) since we are allowed only one integer of a non-zero color and we cannot have $\ell$ consecutive
integers of color 0.  Thus, we see that $\chi$ has form $0^s10^t20^u$ (or its reverse) with $s+t+u=3\ell-4$ and $s,t,u \leq \ell-1$.
As argued previously, we must have $t$ be even, and hence $t = \ell-1$ so that we have $0^s10^{\ell-1}20^t$ with $s+t = 2\ell-3$.
Hence, one of $s$ and $t$ must be $\ell-1$ and the other must be $\ell-2$.  We  assume $s=\ell-1$ (the case $s=\ell-2$
is very similar).
Now that we have $0^{\ell-1}10^{\ell-1}20^{\ell-2}$, consider
$1,4,7,\dots,3\ell-2$ (note that we are using the fact that $\ell \equiv 3 \md{6}$ to end our arithmetic progression
at $3\ell-2$).  Notice that this arithmetic progression consists of integers congruent to $1 \md{3}$ while the
colors $1$ and $2$ are on integers congruent to $0 \md{\ell}$.  Since $\ell \equiv 0 \md{3}$, the
arithmetic progression $1,4,7,\dots,3\ell-2$  is monochromatic of color $0$, and hence is an $\ell$-term
$3$-zero-sum arithmetic progression, a contradiction.  Thus, we can conclude that $\wm(2,\ell;3)=3\ell-2$
for $\ell \equiv 3 \md{6}$.

Lastly, for the $r=3$ case, we consider $\wm(2,\ell;3)$ for $\ell \equiv 1,5 \md{6}$.  For the lower bound, consider the coloring of $[1,3\ell-1]$ given by
$0^{\ell-1}10^{\ell-1}20^{\ell-1}$. As argued above, the only necessary arithmetic progressions to checkare
$\ell$-term ones with common gap 3.  The possibilities are
$1,4,7,\dots,3\ell-2$ and $2,5,8,\dots,3\ell-1$, i.e., the integers congruent to $1$ modulo 3 and the integers congruent
to $2$ modulo 3, respectively.  Since we know that $\ell \equiv 1,5 \md{6}$ we have that one of $\ell$ and $2\ell$ is congruent
to $1$ modulo 3, while the other is congruent to $2$ modulo 3.  Regardless of which is which, we see that neither
of these arithmetic progressions are $3$-zero-sum, thereby proving that $\wm(2,\ell;3) \geq 3\ell$ for
$\ell \equiv 1,5 \md{6}$.

The upper bound is easy in this case since any $3$-coloring of $[1,3\ell]$ using only one 1 and one 2 must have
$3\ell-2$ integers of color $0$, without $\ell$ consecutive integers of color 0.  This is not possible.
Hence, we can conclude that $\wm(2,\ell;3)=3\ell$ for $\ell \equiv 1,5 \md{6}$.

\vskip 5pt
\noindent
{\boldmath $\wm(2,\ell;4)$}. We now move onto $\wm(2,\ell;4)$.  We will start with the lower bounds by giving colorings that avoid
the relevant arithmetic progressions.

For $\ell \equiv 0,2 \md{6}$, we may assume $\ell \geq 6$ since we have calculated $\ell=2$.  
Consider the coloring of $[1,3\ell-1]$ given by $0^{\ell-1}10^{\ell-2}20^{\ell-5}30^4$.  In this situation, $\{\ell, 2\ell-1,3\ell-5\}$, which
are the non-zero colored integers, forms a complete residue system modulo 3.
Hence, any possible $\ell$-term $4$-zero-sum arithmetic progression
$a,a+d,\dots,a+(\ell-1)d$ cannot have $d=3$.  Clearly we do not have
such a progression with $d=1$; so, $d=2$ is the only possibility to check.
Since we require $a+(\ell-1)d \leq 3\ell-1$ we see that $a \leq \ell+1$.
If $a \leq \ell-1$, our progression contains exactly one of
$\ell$ and $2\ell-1$, so that it cannot be $4$-zero-sum.  If $a=\ell$, then
in order for our progression to be $4$-zero-sum, $3\ell-5$ (which has color 3) must be part of the progression.  This means that $a+2j = 3\ell-5$ for
some $j \in \{1,2,\dots,\ell-1\}$.  Since $a=\ell$ this means $2j=2\ell-5$,
which is not possible.  Lastly, if $a=\ell+1$, then our progression has
sum 1 modulo 4 as it contains both $2\ell-1$ and $3\ell-5$.

For $\ell \equiv 3 \md{6}$, we may assume $\ell \geq 9$ since we have calculated $\ell=3$.  
Consider the coloring of $[1,3\ell-1]$ given by $0^{\ell-1}10^{\ell-3}20^{\ell-6}30^6$.  The non-zero colored elements here are
$\ell, 2\ell-2,$ and $3\ell-7$ and these form a complete residue system
modulo 3.  Hence, the only possible $\ell$-term $4$-zero-sum arithmetic
progression $a,a+d,\dots,a+(\ell-1)d$ has $d=2$.  We must still have $a \leq \ell+1$.  If $a \leq \ell$ we use the facts that $\ell$ is odd while
$2\ell-2$ and $3\ell-7$ are even to see that our progression cannot
be 4-zero-sum as every such progression contains at least one of these
elements, but cannot contain both $\ell$ and $3\ell-7$.  If $a=\ell+1$, then
the arithmetic progression contains both $2\ell-2$ and $3\ell-7$ and is
not 4-zero-sum.

For $\ell \equiv 4 \md{6}$, we consider the coloring of $[1,3\ell-1]$ given by $0^{\ell-1}10^{\ell-2}20^{\ell-1}3$.  The non-zero colored elements
here are $\ell$, $2\ell-1$, and $3\ell-1$, the first two being congruent to
1 modulo 3 and the last congruent to 2 modulo 3.  Letting 
$a,a+d,\dots,a+(\ell-1)d$ be an arbitrary $\ell$-term arithmetic progression, consider $d=3$.  If our progression consists of
integers congruent to $1$ modulo 3 then its sum of colors is $1+2\equiv 3 \md{4}$;
if it consists of integers congruent to $2$ modulo 3 then its sum of
colors is $3 \md{4}$; if it consists of integers congruent to $0$ modulo 3,
then $a+(\ell-1)d \geq 3+3(\ell-1)=3\ell>3\ell-1$.  This leaves $d=2$ as
the only possibility.  As above, we have $a \leq \ell+1$.  We also have
than $\ell$ is even while $2\ell-1$ and $3\ell-1$ are odd.  Hence,
for $a \leq \ell$ our progression contains $\ell$ but not $3\ell-1$
or it contains   $2\ell-1$.  In all situations,
our progression is not 4-zero-sum.  If $a=\ell+1$, the progression has
color sum $2+3 \equiv 1 \md{4}$ and, again, the progression is not
$4$-zero-sum.

For $\ell \equiv 1,5 \md{6}$, we consider the coloring of $[1,4\ell-1]$ given by $0^{\ell-1}10^{\ell-1}20^{\ell-1}30^{\ell-1}$.  The non-zero
colored elements $\ell, 2\ell,$ and $3\ell$ form a complete residue
system modulo 3.  Let $a,a+d,\dots,a+(\ell-1)d$ be an arbitrary
$\ell$-term  arithmetic progression. If $d=3$ we have
$a \leq \ell+2$.  If $a \leq \ell$, then our progression contains
exactly one of $\ell$, $2\ell$, and $3\ell$ and is not 4-zero-sum.
If $a = \ell+1$, then $a \equiv 2 \md{3}$ if $\ell \equiv 1 \md{6}$
and $a \equiv 0 \md{3}$ if $\ell \equiv 5 \md{6}$.  In the former case,
our progression contains $2\ell$; in the latter case, our progression
contains $3\ell$.  In either case, we see that the arithmetic progression
is not 4-zero-sum.  If $a = \ell+2$ then we cannot have $\ell \equiv 1 \md{6}$ since then $a \equiv 0 \md{3}$, which tells us that $3\ell$
is part of the progression so that the progression cannot be 4-zero-sum.
Hence, we have $a \equiv 1 \md{3}$ since we have $\ell \equiv 5\md{6}$.
But then $2\ell \equiv 1 \md{3}$ so our progression contains $2\ell$, giving
us that our progression cannot be 4-zero-sum.
Noting that $\ell$ and $3\ell$ are odd, while $2\ell$ is even, consider
$d=2$.  Note that, in this situation, $\ell$ and $3\ell$ cannot both be members of our progression. We must have $a \leq 2\ell+1$ so that our progression clearly contains
at least one of $\ell, 2\ell, 3\ell$ and not both $\ell$ and $3\ell$,
yielding that the arithmetic progression is not 4-zero-sum.
  Hence, if $d=2$ our progression is not 4-zero-sum.
What remains is the case $d=4$.  Here we must have $a \leq 3$.
If $a=1$, or $a=3$, then our progression consists of integers congruent to
$1$ modulo 4, respectively, $3$ modulo 4.  We next note that one of $\ell$ and $3\ell$ is congruent
to $1$ modulo 4 while the other is congruent to 3 modulo 4.  Hence, our
progression cannot be 4-zero-sum.  If $a=2$, then $2\ell$ is a member of
the arithmetic progression so that it is not 4-zero-sum.

We now move onto the upper bounds for $\wm(2,\ell;4)$.

The cases $\ell \equiv 1,5 \md{6}$ are easy so we will
do them first.  Assume, for a contradiction, that there exists
a coloring of $[1,4\ell]$ by $\mathbb{Z}_4$ that does not admit
two terms of the same non-zero color or an $\ell$-term 4-zero-sum
arithmetic progression. Since any such coloring
of $[1,4\ell]$   uses at most one of each non-zero color, it must have
at least $4\ell-3$ integers of color 0. This gives us  
at least $\ell$ consecutive integers of color 0,
a contradiction since these consecutive integers form an
$\ell$-term 4-zero-sum arithmetic progression.

For the cases $\ell \equiv 0,2,3,4 \md{6}$, we will consider the coloring
forms
of $[1, 3\ell]$ given by: (i) $0^s10^t20^u30^v$; (ii) $0^s10^t30^u20^v$;
and (iii) $0^s20^t10^u30^v$, and
leave the reverse colorings' argument details to the reader (which follow
by application of the involution of $[1,3\ell]$ given by $i \mapsto 3\ell+1-i$).  We assume, for a contradiction, that each coloring
avoids the requisite progressions.

For any of the colorings we have: $s,t,u,v \leq \ell-1$ and $s+t+u+v = 3\ell-2$ and that the non-zero colored integers must be a complete residue
system modulo 3 in order to avoid the $\ell$-term 4-zero-sum arithmetic
progressions given by all integers congruent to $i$ modulo 3 for some $i$.

We start with coloring (i): $0^s10^t20^u30^v$.  We know that $s+t+u \geq 2\ell-2$ so that $0^s10^t20^u$ contains the coloring of $[1,2\ell]$.
In order for $1,3,5,\dots,2\ell-1$ and $2,4,\dots,2\ell$ to avoid
being 4-zero-sum, the parity of the   integers colored 1 and 2 must be different.  Similarly, by considering $0^t20^u30^v$ we can deduce
that the parity of the integers colored 2 and 3 must be different.
Hence, the  integers colored 1 and 3 have the same parity.
The integers colored 1 and 3 are $s+1$ and $s+t+u+3$.
If $t+u < 2\ell-3$ then $s+1, s+3,\dots,s+t+u+3,\dots,s+2\ell-1$
is an $\ell$-term 4-zero-sum arithmetic progression.  We can conclude
that $t+u \geq 2\ell-3$.  Since $t,u \leq \ell-1$, we may have
(a) $t=u=\ell-1$; (b) $t=\ell-1$ and $u=\ell-2$; or (c) $t=\ell-2$ and $u=\ell-1$.  

If we have (a), then our non-zero integers are
$s+1, s+\ell+1,$ and $s+2\ell+1$.  Since these must form a complete
residue system modulo 3, we cannot have $\ell \equiv 0,3 \md{6}$. Since
$s+1$ and $s+\ell+1$ must have different parities, we cannot have
$\ell \equiv 2,4\md{6}$.  We conclude that (a) may not occur.

If (b) holds, then $t+1 = \ell$.  Since the non-zero integers
must have different values modulo 3, we cannot have $t+1 \equiv 0 \md{3}$.
Hence, $\ell \not \equiv 0,3 \md{6}$.  Similarly, we cannot
have $u+1 \equiv 0 \md{3}$.  Since $u=\ell-2$ we have $u+1 \equiv \ell-1$
so that $\ell \not \equiv 4 \md{6}$.  Lastly, in order for the integers
colored 1 and 3 to be different modulo 3, we cannot have
$t+u+2 \equiv 0 \md{3}$.  Since $t+u = 2\ell-3$, we cannot have $\ell \equiv 2 \md{6}$.

If (c) holds, essentially the same argument as that for (b) can be employed.

Next, we consider the coloring (ii): $0^s10^t30^u20^v$.
We must have $t =\ell-1$ for otherwise $0^s10^t30^u$ contains
$\ell$ consecutive terms, including the  integers colored
1 and 3; that is, an $\ell$-term 4-zero-sum arithmetic progression.  Further, by considering
$1,3,5,\dots,2\ell-1$ and $2,4,\dots,2\ell$ we see that the
integers colored $1$ and $3$ cannot have the same parity.
This means that $\ell$ must be odd so that $\ell \not \equiv 0,2,4 \md{6}$.
However, if $\ell \equiv 3\md{6}$ we see that these two integers are
equivalent modulo 3, which is also not allowed.

The coloring (iii) can be analyzed by essentially the same argument as that given for (ii).

Having provided matching upper and lower bounds for this last formula, the proof is complete.
\B

We can also provide a formula for the first column of Table 3.  The fact that $\wm(k,2;2)=3$ for all $k$ is trivial.

\begin{thm} Let $k \geq 2$ be an integer. Then
$$
\wm(k,2;3)=\left\{
\begin{array}{ll}
2k -1&\mbox{ if } 2\mid k\\
2k&\mbox{ if }2 \nmid k
\end{array}
\right.
$$
and
$$
\wm(k,2;4)=\left\{
\begin{array}{ll}
3k-2&\mbox{ if } k \equiv 0,2, 3 \md{6}\\
3k-1&\mbox{ if } k \equiv 4 \md{6}\\
3k&\mbox{ if } k \equiv 1, 5 \md{6}
\end{array}
\right.
$$
\end{thm}

\noindent
{\it Proof.}
We first prove the formula for $\wm(k,2;3)$.  For the lower bounds, consider the colorings
$1^{k-1}01^{k-2}$ for $k$ odd and $1^{k-1}01^{k-1}$ for $k$ even.  For the upper bounds,
to avoid $2$-term $3$-zero-sum arithmetic progressions we may only use the color $0$ once
and we cannot have both colors 1 and 2.  Let $c$ be $1$ or $2$.  We cannot have $k$ consecutive integers
of color $c$, so our coloring has form $c^s0c^t$ with $s,t \leq k-1$.  If $k$ is odd, then we
are done since our coloring has maximum length $2k-1$.  If $k$ is even, we cannot have 
$s=t=k-1$ for otherwise the $k$-term arithmetic progressions $1,3,5,\dots,2k-1$ is monochromatic
of color $c$. This completes the proof of the formula
for $\wm(k,2;3)$.

Now onto $\wm(k,2;4)$. For the lower bounds, it is left to the reader to check that the following colorings
avoid the requisite arithmetic progressions:  $1^{k-2}01^{k-1}21^{k-2}$ for $k \equiv 0,2,3 \md{6}$;
$1^{k-1}01^{k-1}21^{k-2}$ for $k \equiv 4 \md{6}$; and
$1^{k-1}01^{k-1}21^{k-1}$ for $k \equiv 1,5 \md{6}$

To finish the proof, we now justify upper bounds for $\wm(k,2;4)$.
We assume, for a contradiction, that in each
case we have a coloring that avoids the relevant arithmetic progressions.
 We can have at most one integer of each of color $0$ and $2$.  Further, we cannot
have both colors $1$ and $3$.  Let $c$ be either $1$ or $3$.  Hence, we can
conclude that any coloring that avoids $k$-term monochromatic arithmetic progressions
of a non-zero color and $2$-term $3$-zero-sum arithmetic progressions has form
$c^s0c^t2c^u$ or its reverse (the situation where we use only one of the colors 0 and 2
is easily dismissed as a possibility since we cannot have length longer than $2k$).  We will only consider this coloring and leave
the reverse coloring's analysis to the reader.

If $k \equiv 1,5 \md{6}$, the argument is essentially identical to the one for $\wm(2,\ell;3)$ in the proof of Theorem \ref{th9} by
replacing $\ell$ with $k$, changing the color $0$ to $c$, and changing the color $1$ to $0$.

If $k \equiv 4 \md{6}$, in order for $c^s0c^t2c^u$ with $s,t,u \leq k-1$ to have length
$3k-1$, we must have $s=t=u=k-1$.
Hence, our coloring is $c^{k-1}0c^{k-1}2c^{k-1}$.  By considering the integers congruent to $2$ modulo 3,
we see that the integers of color $0$ and $2$ are congruent to $0$ modulo $k$, and hence are both
congruent to $1$ modulo 3.  Hence, $2, 5, 7, \dots, 3k-1$ is a $k$-term monochromatic arithmetic progression
of color $c \neq 0$, a contradiction, thereby finishing this case.

If $k \equiv 0,2,3 \md{6}$, in order for $c^s0c^t2c^u$ with $s,t,u \leq k-1$ to have length
$3k-2$, we must have one of $s,t,u$ equal to $k-2$ and the other two equal to $k-1$.
If $s=t=k-1$ and $u=k-2$, then the integers congruent to 1 modulo 3 form a monochromatic
$k$-term arithmetic progression of color $c$ since the integers of color 0 and 2
are congruent to either $0$ or $2$ modulo 3.  If $s=k-2$ and $t=u=k-1$, first consider $k \equiv 0,2 \md{6}$
so that $k$ is even.  Then $k,k+2,k+4,\dots,3k-2$ is a $k$-term monochromatic arithmetic progression
of color $c$.  Next, consider $k \equiv 3 \md{6}$.  In this situation, the integers of color 0 and 2
are both congruent to 2 modulo 3.  Hence, the integers congruent to 1 modulo 3 form a monochromatic
$k$-term arithmetic progression of color $c$.  Lastly, consider $s=u=k-1$ and $t=k-2$.
Then the integers of color 0 and 2 are $k$ and $2k-1$.  Since $k \equiv 3 \md{6}$ we see that
both of these integers are congruent to either $0$ or $2$ modulo 3. Hence, the integers congruent to 1 modulo 3 form a monochromatic
$k$-term arithmetic progression of color $c$.
\B

One final piece we can take from Table 3 concerns $\wm(k,k;k)$.
We have
 $\wm(k,k;k) \geq w(k;2)$.
This holds, since, by definition, there exists a $2$-coloring $\chi: [1,w(k;2)-1] \rightarrow \{0,1\}$ that does not admit
a monochromatic $k$-term arithmetic progression.  Necessarily, we do not have a $k$-term $k$-zero-sum
arithmetic progression since such a progression must be monochromatic.

As was done when considering the $2$-color restriction $\wtwo(k;r)$ of $\w(k;r)$, we investigate
what happens when we restrict
the number of colors to two for these mixed numbers.  This will hopefully allow us to
bound $\wm(k,\ell;r)$.

\begin{defn}
Let $k,\ell,r \in \mathbb{Z}^+$ with $k,\ell,r \geq 2$.  Define $\wmtwo(k,\ell;r)$ to be the minimum integer $n$ such that
any coloring of $[1,n]$ by $\mathbb{Z}_2$ admits either a $k$-term monochromatic arithmetic
progression of  color 1 or an $\ell$-term $r$-zero-sum arithmetic progression.
\end{defn}

Existence is clear since $\wmtwo(k,\ell;r) \leq \wm(k,\ell;r)$ because $\mathbb{Z}_2 \subseteq \mathbb{Z}_r$.

\begin{center}
\begin{tabular}{l|c||c|c|c|c||} 
\multicolumn{1}{c}{}&\multicolumn{1}{c||}{}      &       
\multicolumn{1}{c|}{$\ell=2$}
   &  \multicolumn{1}{c|}{$\ell=3$}  &  \multicolumn{1}{c|}{$\ell=4$}  &
    \multicolumn{1}{c||}{$\ell=5$}  \\ \hline\hline  
\multicolumn{1}{c|}{$k$}&\multicolumn{1}{c||}{$r$}&\multicolumn{1}{c}{}&\multicolumn{1}{c}{}&\multicolumn{1}{c}{}&\\\hline
    & $2$                   &  3  &  6  &  7  &  10      \\ 
   2& $3$                  &3  &6 & 7& 10        \\ 
    &$4$                   & 3 &6 &7 & 10       \\\hline 
    & $2$                  & 3 & 7& 7&  15      \\ 
   3& $3$                  & 6 & 9 &14 & 21       \\ 
    &$4$                   & 6 &9 &18 & 22       \\\hline 
    &$2$                   & 3 & 8& 7&  20      \\ 
   4&$3$                   & 7 & 9& 16 & 23       \\ 
    &$4$                   & 7 &18 &35 &  33      \\\hline 
    &$2$                   & 3 & 8&7 & 21       \\ 
   5&$3$                   &10  & 9&18 & 26       \\ 
    &$4$                   &10  &22 &35 & 37       \\ 
\hline\hline
    \end{tabular}
\vskip 5pt
\centerline{\small {\bf Table 4}:  Values   for $\wmtwo(k,\ell;r)$ for small $k, \ell$, and $r$}
    \end{center}

As we can see, the values here are not as irregular as in Table 3.  However,
as explained in the observations below, attempting to find a formula
or constructive lower bound other than
the apparent $(2k-1)/2k$ and $(2\ell-1)/2\ell$ formula occurring in the first column and first few rows
(which we leave to the reader to investigate), does not seem hopeful.

\vskip 5pt
\noindent
{\bf Observations.}
Via essentially the same argument as that presented for the proof of Proposition \ref{exist}, we have
$\wmtwo(k,k;k) = w(k;2)$.
We also have relationships with the classical van der Waerden numbers in at
least two other ways.  First, we
have $\wmtwo(k,3;t)=w(k,3)$ for all $t \geq 4$, where
$w(k,3)$ is the minimum integer such that any $2$-coloring of $[1,w(k,3)]$ admits either
a monochromatic $k$-term arithmetic progression of the first color or a
monochromatic $3$-term arithmetic progression of the second color. Hence, a result due to Li and  Shu \cite{LS}
gives us that $\wmtwo(k,3;t) > \left( \frac{8}{729}\right) \frac{k^2}{\log^2 k}$ for
sufficiently large $k$ when $t \geq 4$.  Generalizing this, we see that
$\wmtwo(k,\ell;t) = w(k,\ell)$ for all $t \geq \ell+1$.  
Second, we have $\wmtwo(k,\ell;\ell) = \wmtwo(\ell,\ell;\ell)=w(\ell;2)$ for all $k \geq \ell$.
This holds since the first $\ell$ terms of a $k$-term arithmetic progression of color 1 
form an $\ell$-term $\ell$-zero-sum arithmetic progression provided $k \geq \ell$.

\vskip 5pt

Given that many instances of $\wmtwo(k,\ell;r)$ are equal to certain classical van der Waerden numbers,
attempting to
find a formula for these does not seem to be a  good use of time.  

\section{Conclusion and Open Questions} \label{sec4}

\vskip 5pt

If we had a proof of the existence of $\w(k;k)$ (resp., $\wm(k,k;k)$) that did not rely on
the existence of $w(k;k)$, we would have a proof of the existence of $w(k;2)$
by Proposition \ref{exist} (resp., the observation above).  It is an elementary exercise
(see \cite{LR}) to deduce the
existence of $w(k;r)$ from $w(k;2)$ for any $r \in \mathbb{Z}^+$.
Hence, the independent
existence of $\w(k;r)$ (or $\wm(k,k;k)$) implies the existence of $w(k;r)$.
We can state this (for $\w(k;r)$) in the following manner.

\begin{thm} Under the condition that $r \mid k$ the following holds: $\w(k;r)$ exists for all $r$ and $k$   if and only if $w(k;r)$ exists for all $r$ and $k$.
\end{thm}

Unfortunately, all attempts by this author to prove the existence
of $\w(k;r)$ independently from the existence of $w(k;r)$ and its
proofs have  been unsuccessful.

We end with some open questions and problems.

\begin{itemize} \renewcommand{\itemsep}{2pt}

\item[Q1.] Is it true that $\w(k;3)=\wtwo(k;3)$? 

\item[Q2.] Prove or disprove: $\wtwo(k;3)=k^2$.

\item[Q3.] Prove the existence of $\w(k;r)$ and/or $\wm(k,\ell;r)$ independently from van der Waerden's theorem
and its proofs.

\item[Q4.] One useful extension of van der Waerden's theorem is that we can also guarantee that the common gap in
the arithmetic progression has the same color as the arithmetic progression.
Along these lines, when $r \mid k$, investigate
the minimum integer $\wtwo^\ast(k;r)$   such
that every coloring $\chi:[1,\wtwo^\ast(k;r)]\rightarrow \{0,1\}$ admits a $(k-1)$-term arithmetic progression
$a, a+d, a+2d,\dots,a+(k-2)d$ such that $\chi(d) + \sum_{i=0}^{k-2} \chi(a+id) \equiv 0\md r$.
The same can be investigated via an appropriate analogue of $\wm(k,\ell;r)$.

\end{itemize}


\begin{thebibliography}{9} \footnotesize \parskip=1pt

\bibitem{AM}
S. D. Adhikari and E. Mazumdar, 
The polynomial method in the study of zero-sum theorems,
{\it Int. J. Number Theory} {\bf 11} (2015), 1451-1461. 

\bibitem{AC} N. Alon and Y. Caro,
On three zero-sum Ramsey-type problems,
{\it J. Graph Theory} {\bf 17} (1993), 177-192. 

\bibitem{BCRY}
P. Balister, Y. Caro, C. Rousseau, and R. Yuster, 
Zero-sum square matrices, 
{\it European J. Combin.} {\bf 23} (2002), 489-497. 
 
\bibitem{B} A. Bialostocki,  Zero sum trees: a survey of results and open problems, in 
 {\it Finite and Infinite Combinatorics in Sets and Logic}, NATO ASI Series {\bf 411} (Series C: Math. and Physical Sci.), Springer, Dordrecht, 1993.

\bibitem{BBCY} A. Bialostocki, G. Bialostocki, Y. Caro, and R. Yuster,  
Zero-sum ascending waves,
{\it J. Combin. Math. Combin. Comput.} {\bf 32} (2000), 103-114. 

\bibitem{BBS} A. Bialostocki, G. Bialostocki, and D. Schaal, 
A zero-sum theorem,
{\it J. Combin. Theory Ser. A} {\bf 101} (2003), 147-152. 

\bibitem{BD} A. Bialostocki and P. Dierker,
On the Erd\H{o}s-Ginzburg-Ziv theorem and the Ramsey numbers for stars and matchings,
{\it Discrete Math.} {\bf 110} (1992),  1-8. 

\bibitem{BSS} A. Bialostocki, R. Sabar, and D. Schaal, On a zero-sum generalization of a variation of Schur's equation, {\it 
Graphs Combin.} {\bf 24} (2008), 511-518. 

\bibitem{C} Y. Caro,
Zero-sum problems -- a survey,
{\it Discrete Math.} {\bf 152} (1996), 93-113. 
 

\bibitem{EGZ} P. Erd\H{o}s, A. Ginzberg, and A. Ziv,
Theorem in additive number theory,
{\it Bulletin Research Council Israel} {\bf 10F} (1961-2), 41-43.

\bibitem{GG} W. Gao and A. Geroldinger, 
Zero-sum problems in finite abelian groups: a survey,
{\it Expo. Math.} {\bf 24} (2006), 337-369. 



\bibitem{G} D. Grynkiewicz, A weighted Erd\H{o}s-Ginzberg-Ziv theorem,
{\it Combinatorica} {\bf 26} (2006), 445-453.



\bibitem{LR} B. Landman and A. Robertson, \textit{Ramsey Theory on the Integers}, second edition, American Math. Society, 2014.

\bibitem{LS} Y. Li and J. Shu, A lower bound for off-diagonal van der Waerden numbers,
{\it Adv. Appl. Math.} {\bf 44} (2010), 243-247.

\bibitem{R} A. Robertson, Zero-sum generalized Schur numbers, 
preprint at {\tt www.aaronrobertson.org}.
 

\bibitem{vdw} B. L. Van der Waerden, Beweis einer baudetschen Vermutung,
{\it Nieuw Archief voor Wiskunde} {\bf 15} (1927), 212-216.

 





\end{thebibliography}
 \end{document}